\documentclass[12pt]{article}
\usepackage[dvips]{graphicx}
\DeclareGraphicsExtensions{ps,eps}
\makeatletter
\def\section{\@startsection{section}{1}{\z@}{-1.5ex plus -.5ex
minus -2.ex}{1ex plus .2ex}{\large\bf}}
{

\input{amssym.def}
\input{amssym}

\def\rit{{\Bbb R}}

\def\@thmcounterstep{}
\long\def\@makecaption#1#2{\vskip 10pt \setbox\@tempboxa\hbox{#1.#2}
\ifdim \wd\@tempboxa >\hsize
#1.#2\par 
\else99
\hbox to\hsize{\hfil\box\@tempboxa\hfil}
\fi}
\def\ps@headings{
\def\@oddhead{\footnotesize\rm\hfill\runninghead\hfill}
\def\@evenhead{\@oddhead}
\def\@oddfoot{\rm\hfill\thepage\hfill}\def\@evenfoot{\@oddfoot}}
\newtheorem{Theorem}{Theorem}[section]

\newtheorem{Definition}[Theorem]{Definition}
\newtheorem{Proposition}[Theorem]{Proposition}

\makeatother
\title
{
Non relativistic strings may be approximated by relativistic strings
}
\def\runninghead{\quad 
Non relativistic strings may be approximated by relativistic strings
}
\author{
{\em Yann Brenier}\thanks{CNRS, Universit\'e de Nice, FR 2800,
brenier@math.unice.fr}}
\date{} 
\begin{document}
\pagestyle{headings}
\flushbottom
\maketitle
\vspace{-10pt}

\subsection*{}
Keywords : relativistic equations, hyperbolic pdes, 
extremal surfaces, weak completion
\\
\\
AMS classification~: 35Q, 35Q60, 49, 53

\section*{Abstract}

We show that bounded families of global classical relativistic strings
that can be written 
as graphs are relatively compact in $C^0$ topology, but their
accumulation points include many non relativistic strings.
We also provide an alternative formulation of these relativistic
strings and characterize their ``semi-relativistic'' completion.

\section{Some relativistic strings and their non relativistic limits}

Let us consider a graph
$$
(t,s)\in{\bf R}\times{\bf R}\rightarrow (t,s,X(t,s))\in
{\bf R}^2\times{\bf R}^d,
$$
defined by 
a sufficiently smooth (at least locally Lipschitz continuous)
function $X$.
According to string theory (see \cite{Po}, for instance),
this graph defines a global classical relativistic string if and only if,
for all bounded open set $\Omega \subset{\bf R}^2$, 
$X$ makes stationary, with respect
to all perturbations, compactly supported in $\Omega$,
the Nambu-Goto Action defined by
$$
\int_\Omega \sqrt{(1+\partial_s X^2)(1-\partial_t X^2)
+(\partial_t X\cdot\partial_s X)^2}\;\;
dtds
$$
which is nothing but the area (over $\Omega$) of the graph,
in the space ${\bf R}^2\times{\bf R}^d$,
with respect to the Minkowski metric $(-1,+1,\cdot\cdot\cdot,+1)$
(for which the speed of light has unit value).
Since we limit ourself to graphs, we automatically
exclude many kinds of relativistic strings, in particular
loops are ruled out. 
In this limited framework, the 
variational principle just means that $X$ is a solution to
the following
first order partial differential system (of hyperbolic type):
\begin{equation}
\label{string equation}
\partial_t(B\partial_t X-C\partial_s X)
-\partial_s(C\partial_t X+D\partial_s X)=0,
\end{equation}
where
$$ 
B=\frac{1+\partial_s X^2}{A},\;\;\;
C=\frac{\partial_t X\cdot\partial_s X}{A},\;\;\;
D=\frac{1-\partial_t X^2}{A},
$$
$$
A=\sqrt{(1+\partial_s X^2)(1-\partial_t X^2)
+(\partial_t X\cdot\partial_s X)^2}.$$
\\
We say that such a string is global if $X$ is a global solution,
i.e. for the full range $-\infty<t<+\infty$, of (\ref{string equation}).
In the present paper, we exhibit some compactness
properties of these global relativistic strings
and characterize their limits.
\\
In order to motivate this work,
let us first consider,
given some constant $0<\kappa<1$,
the non trivial family ${\bf F_\kappa}$ of
global solutions to the relativistic string equation (\ref{string equation}), 
made of
all $X$ that satisfy, 
the linear wave equation:
\begin{equation}
\label{wave}
\partial_{tt}X=\kappa^2\partial_{ss}X,\;\;\;(t,s)\in {\bf R}^2,
\end{equation}
together with the nonlinear constraint:
\begin{equation}
\label{relativistic condition}
\kappa\partial_{t}X\cdot\partial_{s}X=0,\;\;\;
\partial_{t}X^2+\kappa^2\partial_{s}X^2=1-\kappa^2,
\end{equation}
at time $t=0$.
Any solution of the wave equation (\ref{wave}) which does not satisfy
(\ref{relativistic condition}) will be subsequently
called a non relativistic string.
\\
To check that, indeed, every $X\in{\bf F_\kappa}$ is a global solution
to the relativistic string equation (\ref{string equation}),
let us first notice that
the wave equation (\ref{wave}) also reads:
\begin{equation}
\label{wave bis}
(\partial_{t}\overline{+}\kappa\partial_{s})
(\partial_{t}X\underline{+}\kappa\partial_{s}X)=0,
\end{equation}
which leads to the celebrated d'Alembert formula:
\begin{equation}
\label{dalembert}
(\partial_{t}X\underline{+}\kappa\partial_{s}X)(t,s)
=(\partial_{t}X\underline{+}\kappa\partial_{s}X)(0,s\underline{+}\kappa t).
\end{equation}
Thus, condition (\ref{relativistic condition}), which can be written:
\begin{equation}
\label{relativistic bis}
|\partial_{t}X\underline{+}\kappa\partial_{s}X|^2=1-\kappa^2,
\end{equation}
is propagated by
the wave equation (\ref{wave}) and, therefore, holds true at all
time if it does at time $0$. Finally, we get from
(\ref{relativistic condition}),
$A=\kappa(1+\partial_{s}X^2)$, $B=\kappa^{-1}$, $C=0$ and $D=\kappa$,
so that equation (\ref{string equation}) reduces to (\ref{wave}).
\\
Let us now study ${\bf F_\kappa}$ 
from the viewpoint
of compactness and completeness:

\begin{Theorem}
\label{triv}

The family ${\bf F_\kappa}$ of all $X$ satisfying 
(\ref{wave},\ref{relativistic condition}),
with normalisation $X(0,0)=0$, is a relatively 
compact subset of $C^0({\bf R}^2;\;{\bf R}^d)$.
The closure of ${\bf F_\kappa}$
is made of all functions functions $X$
satisfying (\ref{wave}) and:
\begin{equation}
\label{non relativistic condition}
\partial_{t}X^2+\kappa^2\partial_{s}X^2
+2\kappa|\partial_{t}X\cdot\partial_{s}X|\le 1-\kappa^2.
\end{equation}
\end{Theorem}

The main point of this very easy result (see the proof below)
is that there are many non relativistic strings that
can be uniformly approximated by relativistic strings.
We will call them ``subrelativistic'' strings.
In other words, under completion, algebraic constraints
generated by relativity requirements can
be relaxed as algebraic inequalities.
As we will see below, the situation is very different
for minimal surfaces in Riemannian geometry.
This difference is, unsurprisingly, due to the hyperbolic
character of the string equations, in sharp contrast with the
minimal surface equations, of elliptic nature,
for which elliptic regularity applies.
\\
Before proving Theorem \ref{triv}, let us provide an elementary example
in the case $d=3$. For each integer $n$, we consider the unique solution
$X^{(n)}$ to the linear wave equation (\ref{wave}) with $\kappa=2^{-1/2}$,
and initial conditions:
$$
\partial_t X^{(n)}(0,s)=0,\;\;\;X^{(n)}(0,s)=
$$
$$
(\cos s-1,\frac{\sin(n+1)s}{2(n+1)}+\frac{\sin(n-1)s}{2(n-1)},
\frac{\cos(n+1)s-1}{2(n+1)}+\frac{\cos(n-1)s-1}{2(n-1)}),
$$
This solution satisfies the relativistic constraints 
(\ref{relativistic condition}), since $\kappa^2=1/2$ and
$$
\partial_s X^{(n)}(0,s)=
(-\sin s,\cos s\cos ns,-\cos s\sin ns),
\;\;\;\partial_t X^{(n)}(0,s)=0.
$$
Then, as $n\rightarrow +\infty\;$, $\;X^{(n)}(t,s)$ uniformly converges
toward a limit $X(t,s)$, still solution to
the wave equation (\ref{wave}) with $\kappa=2^{-1/2}$, 
but with initial conditions
$$
\partial_t X(0,s)=0,\;\;\;
X(0,s)=(\cos s-1,0,0),
$$
which makes $(t,s)\rightarrow (t,s,X(t,s))$ a ``subrelativistic'' string,
but not a relativistic one.


\subsection*{Proof of Theorem \ref{triv}}

The proof is elementary. Notice first that
(\ref{relativistic condition}) implies that ${\bf F_\kappa}$ is made of
uniformly Lipschitz functions $X$. With normalization $X(0,0)$=0,
this is enough, according to Ascoli's theorem, to see
that ${\bf F_\kappa}$ is relatively compact for the uniform convergence
on any compact subset of ${\bf R}^2$.
Next, notice that, just as condition
(\ref{relativistic condition}), condition (\ref{non relativistic condition})
can be written
\begin{equation}
\label{non relativistic bis}
|\partial_{t}X\underline{+}\kappa\partial_{s}X|^2\le 1-\kappa^2.
\end{equation}
Thus, both conditions are preserved by the wave equation (\ref{wave bis}).
Let us consider a sequence $X_n$ in ${\bf F_\kappa}$. Up to 
extracting a subsequence, 
we may assume that $X_n(t,s)$ converges to some
limit $X(t,s)$ uniformly on any compact subset of ${\bf R}^2$,
meanwhile $\partial_t X_n$ and $\partial_s X_n$ respectively
converge to $\partial_t X$ and $\partial_s X$ for the
weak-* topology of $L^\infty({\bf R}^2;\;{\bf R}^d)$.
Thus $Y=\partial_s X$ and $W=-\partial_t X$ must take their values
in the closed convex hull of 
\begin{equation}
\label{relativistic manifold}
S=\{(Y,W)\in {\bf R}^{d+d}\;;\;\;\;|W\underline{+}\kappa Y|^2=1-\kappa^2\},
\end{equation}
which exactly is
\begin{equation}
\label{non relativistic manifold}
\{(Y,W)\in {\bf R}^{d+d}\;;\;\;\;|W\underline{+}\kappa Y|^2\le 1-\kappa^2\}.
\end{equation}
Conversely, let us consider a solution $X$ to the wave equation (\ref{wave})
that satisfies $X(0,0)=0$ and (\ref{non relativistic bis}), which
implies
that $(Y_0,W_0)=(\partial_s X,-\partial_t X)(t=0,\cdot)$ is valued
in the closed convex hull of $S$.
Then, at time 0, we may
find a sequence of (smooth) functions $(Y_n^0,W_n^0)$ valued
in $S$ that converges to $(Y_0,W_0)$
for the weak-* topology of $L^\infty({\bf R};\;{\bf R}^d)$.
(This is a well known and very useful property of weak topologies, 
see \cite{Ta}
for instance.)
Let us consider, for each $n$, the unique solution $X_n$ to (\ref{wave}) 
such that
$$
X_n(0,0)=0,\;\;\;(\partial_s X_n,-\partial_t X_n)(t=0,\cdot)=
(Y_n^0,W_n^0).
$$
Then we observe that $X_n$ converges to $X$ uniformly on any compact
subset of ${\bf R}^2$ and satisfies condition (\ref{relativistic condition}).
The proof of Theorem \ref{triv} is now complete.

\subsection*{Comparison with the Euclidean case}

There is no result like Theorem \ref{triv}
in the Riemannian case, with the Euclidean metric $(+1,\cdot\cdot\cdot,+1)$.
In that case, the area of a graph
$$
(t,s)\in \overline\Omega
\rightarrow (t,s,X(t,s))\in
{\bf R}^2\times{\bf R}^d,
$$
where $\Omega$ is a smooth, bounded, connected open subset of ${\bf R}^2$,
is given by:
$$
{\bf A}_\Omega(X)=\int_\Omega \sqrt{(1+\partial_s X^2)(1+\partial_t X^2)
-(\partial_t X\cdot\partial_s X)^2}dsdt.
$$
Then, the minimal surface equation is just:
$$
\partial_t(B\partial_t X+C\partial_s X)
+\partial_s(C\partial_t X+D\partial_s X)=0,
$$
where
$$ 
B=\frac{1+\partial_s X^2}{A},\;\;\;A=\sqrt{(1+\partial_s X^2)(1+\partial_t X^2)
-(\partial_t X\cdot\partial_s X)^2},
$$
$$
C=\frac{\partial_t X\cdot\partial_s X}{A},\;\;\;
D=\frac{1+\partial_t X^2}{A}.
$$
Let us assume that $X$ is harmonic:
$$
\partial_{tt}X+\partial_{ss}X=0,
$$
and $a=\partial_t X\cdot \partial_s X$, $b=\partial_t X^2-\partial_s X^2$
both vanish along $\partial\Omega$. Since $a$ and $b$ are also harmonic,
they must vanish inside $\Omega$. Thus $B=D=1$, $C=0$, and $X$ is
also a solution to the minimal surface equation.
Let us now consider a sequence of such functions $X_n$ and assume that
the restriction of $X_n$ to the boundary $\partial\Omega$
converges to some limit $X_{\partial\Omega}$, say in $C^0(\partial\Omega)$.
Then $X_{\partial\Omega}$ 
has a harmonic extension $X$ and, due to elliptic regularity,
$X_n$ converges to $X$
in $C^\infty(\Omega)$. Thus, $X$ must satisfy
$\partial_t X\cdot \partial_s X=\partial_t X^2-\partial_s X^2=0$,
and, therefore, is still a solution to the minimal surface equation.
So, in this (over)simplified framework,
there is no way to converge to a graph that is not a minimal surface.
Of course, this can be discussed in a much more general framework,
as in \cite{Fe} (chapter 5.4), or, also, in terms of weak continuity
of determinants and polyconvexity (cf. \cite{Ev}, for instance).
To conclude the discussion between riemannian
and lorentzian metrics, 
let us mention reference \cite{GiIs}, where smooth transitions
between riemannian and lorentzian metrics are discussed for 
self-interesting branes in a Lorentzian space time.

\section{The augmented relativistic string equations}

Theorem \ref{triv} is just a motivation
to study more comprehensively
global solutions $X$
to the string equation (\ref{string equation}). Do they
have some compactness properties? What are their limits?
\\
To achieve this goal, we first embed the string equation in a larger,
augmented system.

\begin{Proposition}
\label{derivation}
Let us consider a solution $X$ to the relativistic
string equation (\ref{string equation}).
Then, the following quantities
\begin{equation}
\label{embedding}
\tau=\frac{-L}{1+Y^2},\;\;\;v=\frac{Y\cdot W}{1+Y^2},
\;\;\;\eta=\frac{-L}{1+Y^2}Y,\;\;\;
\zeta=W-\frac{Y\cdot W}{1+Y^2}Y,
\end{equation}
where 
$$
Y=\partial_s X,\;\;\;W=-\partial_t X,\;\;\;
L=-\sqrt{(1+Y^2)(1-W^2)+(Y\cdot W)^2},
$$
are solutions to the ``augmented system'':
$$\partial_t\tau+v\partial_s \tau=\tau\partial_s v,
\;\;\;
\partial_t v+v\partial_s v=\tau\partial_s \tau,
$$
\begin{equation}
\label{non conservative}
\partial_t\eta+v\partial_s \eta=-\tau\partial_s \zeta,
\;\;\;
\partial_t\zeta+v\partial_s \zeta=-\tau\partial_s \eta.
\end{equation}
In addition, they satisfy the following constraints:
\begin{equation}
\label{constraint}
\tau>0,\;\;\;\tau^2+v^2+\eta^2+\zeta^2=1,\;\;\;\tau v=\eta\cdot\zeta.
\end{equation}

\end{Proposition}

Next, we introduce

\begin{Definition}
\label{subrelativistic string}
We call subrelativistic strings all 
solutions $(\tau,v,\eta,\zeta)$ of the augmented
system (\ref{non conservative}) that satisfy the following algebraic
inequalities:
\begin{equation}
\label{subrelativistic condition}
\tau\ge 0,\;\;\tau^2+v^2+\eta^2+\zeta^2+2|\tau v-\eta\cdot\zeta|\le 1.
\end{equation}

\end{Definition}

We will say that a subrelativistic string $(\tau,v,\eta,\zeta)$
is global whenever it is a global solution to the augmented system,
i.e. for $-\infty<t<+\infty$.
We will see later that a 
necessary condition for $(\tau,v,\eta,\zeta)$
to be a global subrelativistic string 
is the existence of a real constant $\alpha$ such that
$$
\tau\underline{+}(v-\alpha)>0
$$
holds true at $t=0$,
meanwhile a sufficient condition is the further existence of some
constant $\delta>0$ such that:
$$
\delta\le \tau\underline{+}(v-\alpha)\le \frac{1}{\delta}.
$$
Our main result, which will be precisely stated as Theorem \ref{compactness},
asserts, roughly speaking, that global subrelativistic strings form
a natural completion for global relativistic strings.
The main steps of the analysis will be:
\\
1) an almost explicit resolution of the augmented system
for a large class of initial data, thanks to d'Alembert formula,
\\
2) a weak convergence argument, using that (\ref{subrelativistic condition})
defines the closed convex hull of (\ref{constraint}).
\\
Let us finally mention, before proving Proposition \ref{derivation},
that it has been known for a long time that relativistic
string equations can be solved using d'Alembert's formula. (Just like
minimal surfaces can be reduced to harmonic functions.) See 
\cite{Po}, for instance. It turns out that this is also true for 
generalized strings described by the augmented system.

\subsection*{Proof of Proposition \ref{derivation}}

Let
\begin{equation}
\label{Y-W}
Y=\partial_s X,\;\;\;W=-\partial_t X.
\end{equation}
Then the relativistic string equation may be equivalently 
obtained by varying the
Action $\int L(Y,W)dtds$ where $Y$ and $W$ are subject to
$$
\partial_t Y+\partial_s W=0,
$$
and the Lagrangian density $L$ is given by
$$
L(Y,W)=-\sqrt{(1+Y^2)(1-W^2)+(W\cdot Y)^2}.
$$
The resulting equations are
$$
\partial_t Y+\partial_s W=0,
\;\;\;
\partial_t Z+\partial_s V=0,
$$
where $Z$ and $V$ are defined by
$$
Z=\frac{\partial L}{\partial W}(Y,W)=\frac{(1+Y^2)W-(W\cdot Y)Y}{-L}.
$$
$$
V=-\frac{\partial L}{\partial Y}(Y,W)=\frac{(1-W^2)Y-(W\cdot Y)W}{-L}.
$$
In order to write $W$ and $V$ as functions of the evolution variables
$Y$ and $Z$, we introduce 
the hamiltonian function $h$ defined as the partial Legendre transform
$$
h(Y,Z)=\sup_{W\in\rit^d} Z\cdot W-L(Y,W)=\sqrt{1+Y^2+Z^2+(Y\cdot Z)^2}.
$$
Let us introduce
\begin{equation}
\label{q}
q=Y\cdot Z.
\end{equation}
Thus
\begin{equation}
\label{h}
h=\sqrt{1+Y^2+Z^2+q^2},
\end{equation}
\begin{equation}
\label{V-W}
V=\frac{\partial h}{\partial Y}(Y,Z)=\frac{Y+qZ}{h},
\;\;\;\;
W=\frac{\partial h}{\partial Z}(Y,Z)=\frac{Z+qY}{h}.
\end{equation}
The relativistic string equation now reads:
\begin{equation}
\label{relativistic}
\partial_t Y+\partial_s(\frac{Z+qY}{h})=0,
\;\;\;
\partial_t Z+\partial_s(\frac{Y+qZ}{h})=0, 
\end{equation}
where $q$ and $h$ are defined by (\ref{q},\ref{h}).

Next, we follow an idea used in \cite{Br} for the Born-Infeld system
(for which we also refer to \cite{BDLL,Gi,Se,Se2,Se3}),
by adding  to system (\ref{relativistic}) two additional conservation laws, 
for $h$ and $q$, respectively:
\begin{equation}
\label{conservation h}
\partial_t h+\partial_s q=0,
\end{equation}
\begin{equation}
\label{conservation q}
\partial_t q+\partial_s (\frac{q^2-1}{h})=0.
\end{equation}
System (\ref{conservation h},\ref{conservation q}) is known under many
different names, such as the Chaplygin gas equation, the (one-dimensional)
Born-Infeld equations or the Eulerian version of the linear wave equation.
\cite{BDLL}, \cite{Se}, \cite{Se2}. As we will see in the
next section,
this system can be easily integrated by using d'Alembert's formula.
Let us now establish equations (\ref{conservation h},\ref{conservation q})
from the string equation written in form
(\ref{q},\ref{h},\ref{V-W},\ref{relativistic}).
We first get
$$
\partial_t h=\frac{\partial h}{\partial Y}\cdot\partial_t Y
+\frac{\partial h}{\partial Z}\cdot\partial_t Z
$$
$$
=-V\cdot\partial_s W-W\cdot\partial_s V=-\partial_s(W\cdot V)
$$
where
$$
W\cdot V=\frac{(Z+(Y\cdot Z)Y)(Y+(Y\cdot Z)Z)}{h^2}=Z\cdot Y,
$$
which leads to (\ref{conservation h}).
Next, we have
$$
-\partial_t q=Z\cdot\partial_s W+Y\cdot\partial_s V
=Z\cdot\partial_s(\frac{Z+(Y\cdot Z)Y}{h})
+Y\cdot\partial_s(\frac{Y+(Y\cdot Z)Z}{h}).
$$
$$
=Z\cdot\partial_s(\frac{(Y\cdot Z)Y}{h})
+Y\cdot\partial_s(\frac{(Y\cdot Z)Z}{h})
+h\partial_s{\frac{Z^2+Y^2}{2h^2}}.
$$
Observe that
$$
\partial_s\frac{q^2}{h}
=\partial_s\frac{(Y\cdot Z)^2}{h}
=Z\cdot\partial_s\frac{Y(Y\cdot Z)}{h}
+\partial_s Z\cdot\frac{Y(Y\cdot Z)}{h}
$$
$$
=Z\cdot\partial_s\frac{Y(Y\cdot Z)}{h}
+Y\cdot\partial_s\frac{Z(Y\cdot Z)}{h}
-Y\cdot Z\partial_s\frac{Y\cdot Z}{h}
$$
$$
=Z\cdot\partial_s\frac{Y(Y\cdot Z)}{h}
+Y\cdot\partial_s\frac{Z(Y\cdot Z)}{h}
-h\partial_s\frac{(Y\cdot Z)^2}{2h^2}.
$$
Thus
$$
\partial_t q+\partial_s \frac{q^2}{h}=
-h\partial_s\frac{Z^2+Y^2+(Y\cdot Z)^2}{2h^2}
=-h\partial_s\frac{h^2-1}{2h^2}
=\partial_s\frac{1}{h},
$$
which is just (\ref{conservation q}).
\\
Let us finally introduce the rescaled variables:
\begin{equation}
\label{rescaled}
\tau=\frac{1}{h},\;\;\;v=\frac{q}{h},\;\;\;\eta=\frac{Y}{h},\;\;\;
\zeta=\frac{Z}{h}.
\end{equation}
Because of (\ref{q},\ref{h}), they must satisfy 
$$
\tau>0,\;\;\;\tau^2+v^2+\eta^2+\zeta^2=1,\;\;\;\tau v=\eta\cdot\zeta,
$$
which exactly is (\ref{constraint}).
\\
After straightforward calculations,
the ``augmented'' system 
(\ref{relativistic},\ref{conservation h},\ref{conservation q})
can be written in terms of $\tau,v,\eta,\zeta$:
$$
\partial_t\tau+v\partial_s \tau=\tau\partial_s v,
\;\;\;
\partial_t v+v\partial_s v=\tau\partial_s \tau,
$$
$$
\partial_t\eta+v\partial_s \eta=-\tau\partial_s \zeta,
\;\;\;
\partial_t\zeta+v\partial_s \zeta=-\tau\partial_s \eta,
$$
which is nothing but (\ref{non conservative}).
Thus, the proof of Proposition \ref{derivation} is now complete.

\subsection*{Comments on the augmented system}

The augmented system (\ref{non conservative})
makes sense for all $U=(\tau,v,\eta,\zeta)\in{\bf R}^{1+1+d+d}$,
even if (\ref{constraint}) is not satisfied. 
(Notice that $\tau$ may even change sign!)
As a matter of fact system (\ref{non conservative})
can be written as 
$$
\partial_t U+A(U)\partial_s U=0,
$$
where $A(U)$ is a symmetric matrix.
Therefore, this system
is a symmetric hyperbolic system of first order PDEs.
As a consequence,
the Cauchy problem, with initial data at time $t=0$,
is solvable in a neighborhood of $t=0$, for all smooth initial
data $s\in{\bf R}\rightarrow U_0(s)\in {\bf R}^{1+1+d+d}$,
with appropriate behaviour near $s=\underline{+}\infty$.
\\
Surprisingly enough, the augmented system (\ref{non conservative})
is Galilean invariant, under the following transform:
\begin{equation}
\label{galilean}
(t,s)\rightarrow (t,s+ut),\;\;\;
(\tau,v,\eta,\zeta)\rightarrow (\tau,v+u,\eta,\zeta),
\end{equation}
where $u\in {\bf R}$ is a fixed velocity.
Observe that this transform, which is certainly 
ruled out by the relativistic constraint (\ref{constraint}),
is compatible with the ``subrelativistic'' condition 
(\ref{subrelativistic condition}),
provided $|u|$ is not too large.

\subsection*{Comment on relativistic and non-relativistic strings}

As shown in Proposition \ref{derivation},
we can attach a solution $U=(\tau,v,\eta,\zeta)$ of the augmented system
(\ref{non conservative}) to each graph
$(t,s)\rightarrow (t,s,X(t,s))$ corresponding
to a relativistic string, through (\ref{embedding}).
Of course, by construction of the augmented system,
such solutions automatically satisfy constraint (\ref{constraint}).
Conversely, given a smooth solution $U=(\tau,v,\eta,\zeta)$ to
the augmented system (\ref{non conservative}), such that $\tau>0$,
we may define $X(t,s)$ (up to a normalization such as $X(0,0)=0$) by
$$
\partial_s X=\frac{\eta}{\tau},\;\;\;\partial_t X=-\zeta-v\frac{\eta}{\tau}.
$$
Then, if $U$ satisfies (\ref{constraint}), we can check from
(\ref{non conservative}) that, indeed, $X$ solves the
relativistic string equation (\ref{string equation}).

Other solutions to the augmented system (\ref{non conservative}) 
may describe graphs 
$(t,s)\rightarrow (t,s,X(t,s))$ that
are not necessarily relativistic strings.
For instance, consider a non relativistic string, for which
$X$ solves the wave equation (\ref{wave}) but not
necessarily equation (\ref{string equation}).
Then, assuming
$v=0$, $\tau=\kappa$, 
in system (\ref{non conservative}), we get
$$
\partial_t\eta+\kappa\partial_s \zeta=
\partial_t\zeta+\kappa\partial_s \eta=0,
$$
and, by setting:
$$
\eta=\kappa\partial_s X,\;\;\;\zeta=-\partial_t X,
$$
we recover the wave equation (\ref{wave}).
Such a string is relativistic only if (\ref{constraint}) is
satisfied, which means
$$
\kappa^2+\eta^2+\zeta^2=1,\;\;\;\eta\cdot\zeta=0,
$$
or, in other words,
$$
\kappa^2|\partial_s X|^2+|\partial_t X|^2=1-\kappa^2,\;\;\;
\kappa\partial_s X\cdot\partial_t X=0,
$$
which exactly is condition (\ref{relativistic condition}).

\section{Integration of the augmented system}

The augmented system (\ref{non conservative}) 
can also be written in ``diagonal'' form:
\begin{equation}
\label{diagonal}
D_t^\epsilon(v-\epsilon\tau)=0,
\;\;\;D_t^\epsilon(\eta+\epsilon\zeta)=0,
\end{equation}
where $\epsilon\in \{-1,+1\}$ and
\begin{equation}
\label{D-epsilon}
D_t^\epsilon=\partial_t+(v+\epsilon\tau)\partial_s.
\end{equation}
It follows that, $\epsilon$ being fixed in $\{-1,+1\}$,
for any real function $f$ and any constant $r$,
the level set
$$
\{U=(\tau,b,\eta,\zeta);\;\;\;f(v-\epsilon\tau,\eta+\epsilon\zeta)=r\},
$$
is an invariant set for system (\ref{diagonal}).
As a consequence, the following sets are also invariant:
\begin{equation}
\label{well posed manifold}
G_{\alpha,\delta}=\{U=(\tau,v,\eta,\zeta)\in {\bf R}^{1+1+d+d};\;\;\;
\delta\le \tau\underline{+}(v-\alpha)\le \frac{1}{\delta}\;\},
\end{equation}
for any constants $\alpha\in{\bf R}$ and $0<\delta<1$.
Observe that $G_{\alpha,\delta}$ is included in 
$$
\{U=(\tau,v,\eta,\zeta);\;\;\;
\delta\le \tau\le \frac{1}{\delta}\;\}.
$$
Other invariant sets are:
$$
M_\epsilon=\{(\tau,v,\eta,\zeta)\in {\bf R}^{1+1+d+d}\;;
\;\;\;(v+\epsilon \tau)^2+|\eta-\epsilon\zeta|^2=1\}
$$
for $\epsilon=-1,+1$,
as well as their intersection:
\begin{equation}
\label{bi manifold}
M=\{(\tau,v,\eta,\zeta)\in {\bf R}^{1+1+d+d};\;\;\;
\tau v=\eta\cdot\zeta,\;\;\;\tau^2+v^2+\eta^2+\zeta^2=1\},
\end{equation}
which precisely corresponds to the relativistic string constraint 
(\ref{constraint}).
\\
It is now easy to integrate system (\ref{diagonal}) for solutions
valued in $G_{\alpha,\delta}$.

\begin{Proposition}

Let $\alpha\in {\bf R}$, $\delta>0$ be fixed constants.
All solutions $U=(\tau,v,\eta,\zeta)$ to the augmented system
(\ref{non conservative}) valued in the invariant set
$G_{\alpha,\delta}$ (\ref{well posed manifold}) are global and implicitly
defined by:
$$
(v\underline{+}\tau)(t,\xi(t,y))=
(v\underline{+}\tau)(0,\xi(0,y\overline{+} t))),\;\;\;\forall(t,y)\in{\bf R}^2
$$
\begin{equation}
\label{solution}
(\eta\overline{+}\zeta)(t,\xi(t,y))=
(\eta\overline{+}\zeta)(0,\xi(0,y\overline{+} t)),
\end{equation}
where, for each t, $y\rightarrow \xi(t,y)$ is 
a bi-Lipschitz homeomorphism of the real line, with
\begin{equation}
\label{definition xi}
\partial_y\xi(t,y)=\tau(t,\xi(t,y)),
\end{equation}
valued in $[\delta,\frac{1}{\delta}]$.
In addition, $\xi$ is completely determined by $\xi(0,0)=0$ and:
\begin{equation}
\label{solution xi}
(\partial_t\xi\underline{+}\partial_s\xi)(t,y)
=v(0,\xi(0,y\underline{+} t))\underline{+}\tau(0,\xi(0,y\underline{+} t)).
\end{equation}

\end{Proposition}

\subsection*{Proof}

Since $0<\delta\le \tau\le \delta^{-1}$, 
we can define  $\xi(t,y)$
for all $(t,y)\in {\bf R}^2$ by (\ref{definition xi})
in such a way that, in addition,
$$
\partial_t\xi(t,y)=v(t,\xi(t,y)),\;\;\;
\partial_{tt}\xi=\partial_{yy}\xi,
$$
hold true. This is possible, due to the 
two first equations of the augmented system 
(\ref{non conservative}).
Of course, we can normalize $\xi(0,0)=0$.
Next,
$$
(\partial_{t}\overline{+}\partial_s)
(\partial_{t}\xi\underline{+}\partial_s\xi)
$$
follows, and we deduce (\ref{solution xi}) from d'Alembert's formula.
Notice that (\ref{solution xi}) and $\xi(0,0)=0$ entirely 
determine $\xi$, given $\tau$ and $v$ at time $t=0$.
\\
Then, using (\ref{diagonal}), we get:
$$
\partial_t[(v\underline{+}\tau)(t,\xi(t,y\underline{+} t))]=0,
\;\;\;\partial_t[(\eta\overline{+}\zeta)(t,\xi(t,y\underline{+} t))]=0,
$$
which leads to formula (\ref{solution}) and completes the proof.

\subsection*{Generalized solutions}

Just as d'Alembert's formula does for the linear wave equation,
formulae (\ref{solution xi},\ref{solution}) 
provide a natural notion of (global) generalized solutions
for the augmented system (\ref{non conservative}), 
globally and uniquely defined
for $each$ Lebesgue measurable 
initial condition valued in $G_{\alpha,\delta}$, which means
\begin{equation}
\label{initial}
\delta\le \tau(0,\cdot)\underline{+}(v(0,\cdot)-\alpha)\le \frac{1}{\delta}\;,
\end{equation}
for some constants $\alpha\in{\bf R}$, $0<\delta<1$.
Of course, 
for each smooth initial condition, the corresponding generalized
solution automatically is a classical, global, smooth solution to
the augmented system (\ref{non conservative}).
As a matter of fact, condition (\ref{initial}) is nearly a necessary
condition to define a global solution $(\tau,v,\eta,\zeta)$ to
system (\ref{non conservative}) with $\tau>0$.
Indeed, because of (\ref{definition xi}),
$\;\partial_y\xi(t,y)$ must stay positive for all $(t,y)\in{\bf R}^2$. 
Because of d'Alembert formula (\ref{solution xi}), this is possible 
only if
$$
v(0,s)-\tau(0,s)<v(0,s')+\tau(0,s'),\;\;\;\forall s,s'\in {\bf R},
$$
which exactly means 
$\tau(0,\cdot)\underline{+}(v(0,\cdot)-\alpha)>0$,
for some constant $\alpha\in {\bf R}$.
In the rest of the paper, we will consider only generalized solutions
valued in one of the $G_{\alpha,\delta}$.

\subsection*{Weak form}

For each generalized solution valued in $G_{\alpha,\delta}$,
$y\rightarrow \xi(t,y)$ is a bi-Lipschitz homeomorphism of ${\bf R}$, 
for each $t\in {\bf R}$, since
$0<\delta\le \partial_y\xi\le \delta^{-1}<+\infty$. 
The inverse of $\xi(t,\cdot)$ is denoted
by $\xi^{-1}(t,\cdot)$. This allows us to write
(\ref{definition xi})
in the following ``weak'' form:
\begin{equation}
\label{w-definition xi}
\int_{-\infty}^{+\infty} \frac{g(s)}{\tau(t,s)}ds=
\int_{-\infty}^{+\infty} g(\xi(t,y))dy,
\end{equation}
for all functions $g\in L^1({\bf R})$.
Similarly, (\ref{solution xi},\ref{solution}) reads:
\begin{equation}
\label{w-solution xi}
\int_{-\infty}^{+\infty}
(\partial_t\xi\underline{+}\partial_s\xi)(t,y)g(y)dy
=\int_{-\infty}^{+\infty}
\frac{v\underline{+}\tau}{\tau}(0,s)
g(\xi^{-1}(0,s)\overline{+} t)ds,
\end{equation}
$$
\int_{-\infty}^{+\infty}
\frac{v\underline{+}\tau}{\tau}(t,s)g(s)ds=
\int_{-\infty}^{+\infty}
\frac{v\underline{+}\tau}{\tau}(0,s)
g(\xi(t,\xi^{-1}(0,s)\overline{+} t))ds,
$$
\begin{equation}
\label{w-solution}
\int_{-\infty}^{+\infty}
\frac{\eta\overline{+}\zeta}{\tau}(t,s)g(s)ds=
\int_{-\infty}^{+\infty}
\frac{\eta\overline{+}\zeta}{\tau}(0,s)
g(\xi(t,\xi^{-1}(0,s)\overline{+} t))ds,
\end{equation}
for all functions $g\in L^1({\bf R})$.
Using the original variables $(h,q,Y,Z)$ instead of $(\tau,v,\eta,\zeta)$,
we respectively get:
\begin{equation}
\label{ww-definition xi}
\int_{-\infty}^{+\infty} {g(s)}{h(t,s)}ds=
\int_{-\infty}^{+\infty} g(\xi(t,y))dy,
\end{equation}
\begin{equation}
\label{ww-solution xi}
\int_{-\infty}^{+\infty}
(\partial_t\xi\underline{+}\partial_s\xi)(t,y)g(y)dy
=\int_{-\infty}^{+\infty}
(q(0,s)\underline{+}1)
g(\xi^{-1}(0,s)\overline{+} t)ds,
\end{equation}
$$
\int_{-\infty}^{+\infty}
(q(t,s)\underline{+}1)g(s)ds=
\int_{-\infty}^{+\infty}
(q(0,s)\underline{+}1)
g(\xi(t,\xi^{-1}(0,s)\overline{+} t))ds,
$$
\begin{equation}
\label{ww-solution}
\int_{-\infty}^{+\infty}
(Y\overline{+} Z)(t,s)g(s)ds=
\int_{-\infty}^{+\infty}
(Y\overline{+} Z)(0,s)
g(\xi(t,\xi^{-1}(0,s)\overline{+} t))ds,
\end{equation}
for all functions $g\in L^1({\bf R})$.

\section{Weak completion of global relativistic strings}

In this last section, we study the subset
of all global generalized solutions to the augmented system 
(\ref{non conservative}) valued in the invariant subset $G_{\alpha,\delta}$
(defined by (\ref{well posed manifold})
for some fixed constants $\alpha\in{\bf R}$, $0<\delta<1$,
which, in addition, satisfy the relativistic constraint (\ref{constraint}),
or, in other words, are valued in the invariant region $M$
defined by (\ref{bi manifold}), and, therefore, correspond to global
relativistic strings.

We call $\Sigma_{\alpha,\delta}$ the set of
all such solutions.
We also denote:
$$
M_{\alpha,\delta}=M\;\cap\;G_{\alpha,\delta},
$$
i.e.
$$
M_{\alpha,\delta}
=\{U=(\tau,v,\eta,\zeta);
\;\;\delta\le \tau\underline{+}(v-\alpha)\le \frac{1}{\delta}\;;
$$
\begin{equation}
\label{set}
\tau v=\eta\cdot\zeta,\;\;\tau^2+v^2+\eta^2+\zeta^2=1\}.
\end{equation}
An equivalent definition is:
\begin{equation}
\label{set bis}
M_{\alpha,\delta}=
\{\delta\le \tau\underline{+}(v-\alpha)\le \frac{1}{\delta}\;,
\;\;\;(v\underline{+}\tau)^2+|\eta\overline{+}\zeta|^2= 1\;\}\;.
\end{equation}

From the topological point of view, we confer to
$\Sigma_{\alpha,\delta}$ the topology induced by the space
$C^0({\bf R};L^{\infty}_{weak^*}({\bf R};{\bf R}^{1+1+d+d}))$
through the one-to-one transform 
\begin{equation}
\label{transform}
T\;:U=(\tau,v,\eta,\zeta)
\rightarrow 
u=(h,q,Y,Z)=\frac{1}{\tau}(1,v,\eta,\zeta),
\end{equation}
defined on ${\bf R_+}\times{\bf R}^{1+d+d}$.
More precisely, we say that $U_n=(\tau_n,v_n,\eta_n,\zeta_n)$
converges to  $U=(\tau,v,\eta,\zeta)$ if and only if
$TU_n$ converges to $TU$ in
$$
C^0({\bf R};L^{\infty}_{weak^*}({\bf R};{\bf R}^{1+1+d+d})),
$$
i.e.
\begin{equation}
\label{convergence}
\int_{-\infty}^{+\infty} (h_n-h,q_n-q,Y_n-Y,Z_n-Z)(s)g(s)ds\rightarrow 0
\end{equation}
uniformly in $t$ on any compact subset of ${\bf R}$,
for all functions $g\in L^1({\bf R})$,
or, equivalently
\begin{equation}
\label{convergence bis}
\int_{-\infty}^{+\infty} \{\frac{(1,v_n,\eta_n,\zeta_n)}{\tau_n}
-\frac{(1,v,\eta,\zeta)}{\tau}\}(s)g(s)ds\rightarrow 0.
\end{equation}
Notice that $T$ and its inverse
$$
\;\;\;\;T^{-1}\;:u=(h,q,Y,Z)\rightarrow U=(\tau,v,\eta,\zeta)=\frac{1}{h}(h,q,Y,Z)
$$
(which was already used for definition (\ref{rescaled})),
both 
preserve straight lines and convexity on ${\bf R_+}\times {\bf R}^{1+d+d}$.

\begin{Theorem}
\label{compactness}
The set $\Sigma_{\alpha,\delta}$ is relatively compact for the
toplogy of 
$$C^0({\bf R};L^{\infty}_{weak^*}({\bf R};{\bf R}^{1+1+d+d})),$$
induced by $T$ (defined by (\ref{transform},\ref{convergence})).
Its closure 
is the set of all generalized solutions to the augmented system 
(\ref{non conservative}),
in the sense of (\ref{solution xi},\ref{solution}),
valued in $CM\cap G_{\alpha,\delta}$, where
\begin{equation}
\label{abi manifold}
CM=\{U=(\tau,v,\eta,\zeta)\in {\bf R}^{1+1+d+d};\;\;
\tau^2+v^2+\eta^2+\zeta^2+2|\tau v-\eta\cdot\zeta|\le 1\}.
\end{equation}

\end{Theorem}

\subsubsection*{Comment on Theorem \ref{triv}}

Theorem (\ref{compactness}) has Theorem (\ref{triv})
as a corollary. Indeed, let us consider
a solution $X$ to the wave equation (\ref{wave}),
with $0<\kappa<1$.
Assume that $X$ satisfies (\ref{non relativistic condition})
and define $\tau=\kappa$, $v=0$,
$$
\eta=\kappa\partial_s X,\;\;\;\zeta=-\partial_t X.
$$
Then $U=(\tau,v,\eta,\zeta)$ is valued in $CM\cap G_{\alpha,\delta}$, 
for $\alpha=0$ and $\delta=\kappa>0$.
Thus, $U$ can be approximated by a sequence $U_n$
valued in $M\cap G_{\alpha,\delta}$, which means that there is
a sequence of relativistic strings $(t,s)\rightarrow X_n(t,s)$,
such that 
$$
\frac{(1,v_n,\eta_n,\zeta_n)}{\tau_n}\rightarrow
\frac{(1,v,\eta,\zeta)}{\tau}
$$
in $C^0({\bf R};L^{\infty}_{weak^*}({\bf R};{\bf R}^{1+1+d+d}))$,
which, in particular, implies that $X_n(t,s)$ converges to $X(t,s)$
uniformly on all compact subset of ${\bf R^2}$.

\subsubsection*{Proof}

Our proof is elementary and based on closed formulae
(\ref{ww-definition xi},\ref{ww-solution}).
Alternative proofs, based on the Murat-Tartar ``div-curl'' lemma
\cite{Ta}, are possible, following Serre's analysis of the one-dimensional
Born-Infeld equation
\cite{Se}.
\\
Let us consider a sequence $U_n=(\tau_n,v_n,\eta_n,\zeta_n)$
in $\Sigma_{\alpha,\delta}$
and the corresponding variables $TU_n=(h_n,q_n,Y_n,Z_n)$.
Using definitions (\ref{well posed manifold},\ref{bi manifold})
and formulae (\ref{solution}), we have:
$$
\delta\le \tau_n\underline{+}(v_n-\alpha)\le \frac{1}{\delta}\;;
$$
$$
\tau^2_n+v^2_n+\eta_n^2+\zeta_n^2=1,\;\;\;\delta\le \tau_n\le 1,
$$
$$
(v_nn\underline{+}\tau_n)(t,\xi_n(t,y\underline{+} t))=(v_n\underline{+}\tau_n)
(0,\xi_n(0,y))),
$$
$$
(\eta\overline{+}\zeta_n)(t,\xi_n(t,y\underline{+} t))=
(\eta\underline{+}\zeta_n)(0,\xi_n(0,y)),
$$
where
$$
\partial_s\xi_n(t,y)=\tau_n(t,\xi_n(t,y)),
\;\;\;\partial_t\xi_n(t,y)=v_n(t,\xi_n(t,y)),
$$
with normalization $\xi_n(0,0)=0$.
So, we immediately get:
$$
(\partial_t\xi_n)^2+(\partial_s\xi_n)^2\le 1,\;\;\;\partial_s\xi_n\ge \delta.
$$
and deduce, using Ascoli's theorem, that
$\xi_n$ is relatively compact in  $C^0({\bf R}^2)$.
Thus, up to the extraction of a subsequence, $\xi_n(t,y)$ 
uniformly converges to a limit
$\xi(t,y)$ on any compact subset of ${\bf R}^2$.
Since $\delta\le \partial_y\xi_n\le \delta^{-1}$,
we also have $\xi_n^{-1}(t,s)\rightarrow\xi^{-1}(t,s)$ 
uniformly in $(t,s)$, on any compact subset of ${\bf R}^2$.
\\
Let us now consider the initial values $U_n(0,\cdot)$.
Because of constraint (\ref{constraint}), these functions are bounded in 
sup norm. Thus, the sequence $TU_n(0,\cdot)$, 
where $T$ is defined by (\ref{transform}), is bounded in
the space $L^\infty({\bf R};{\bf R}^{1+1+d+d})$. So, up to the extraction
of a further subsequence, we may assume that they converge (in the weak-*
sense)
to some limit $u_{in}=(h_{in},q_{in},Y_{in},Z_{in})$.
Since $T$ and $T^{-1}$ preserve convexity,
$U_{in}=T^{-1}u_{in}$
is valued in the closed convex hull of $M_{\alpha,\delta}$,
that we denote by $CM_{\alpha,\delta}$.
Since $G_{\alpha,\delta}$ is a closed, convex subset
of ${\bf R}^{1+1+d+d}$ and contains
$M_{\alpha,\delta}$, according to definitions
(\ref{well posed manifold},\ref{set}), it certainly contains
$CM_{\alpha,\delta}$.
Thus $U_{in}$ is valued in $CM_{\alpha,\delta}$ and $G_{\alpha,\delta}$.
\\
Let us go back to $TU_n=(h_n,q_n,Y_n,Z_n)$. 
Because of (\ref{ww-solution xi}), we 
have
$$
\int_{-\infty}^{+\infty}
(\partial_t\xi_n\underline{+}\partial_s\xi_n)(t,y)g(y)dy
=\int_{-\infty}^{+\infty}
(q_n(0,s)\underline{+}1)
g(\xi_n^{-1}(0,s)\overline{+} t)ds,
$$
for all functions $g\in L^1({\bf R})$ and $t\in{\bf R}$.
We deduce,
after letting $n\rightarrow \infty$,
$$
\int_{-\infty}^{+\infty}
(\partial_t\xi\underline{+}\partial_s\xi)(t,y)g(y)dy
=\int_{-\infty}^{+\infty}
(q_{in}(s)\underline{+}1)
g(\xi^{-1}(0,s)\overline{+} t)ds.
$$
Next, we use
that $TU_n$ satisfies
(\ref{ww-definition xi},\ref{ww-solution}):
$$
\int_{-\infty}^{+\infty} {g(s)}{h_n(t,s)}ds=
\int_{-\infty}^{+\infty} g(\xi_n(t,y))dy,
$$
$$
\int_{-\infty}^{+\infty}
(q_n\underline{+}1)(t,s)g(s)ds=
\int_{-\infty}^{+\infty}
(q_n\underline{+}1)(0,s)
g(\xi_n(t,\xi_n^{-1}(0,s)\overline{+} t))ds,
$$
$$
\int_{-\infty}^{+\infty}
(Y_n\overline{+} Z_n)(t,s)g(s)ds=
\int_{-\infty}^{+\infty}
(Y_n\overline{+} Z_n)(0,s)
g(\xi_n(t,\xi_n^{-1}(0,s)\overline{+} t))ds,
$$
for all functions $g\in L^1({\bf R})$ and $t\in{\bf R}$.
As $n\rightarrow +\infty$, each right-hand side 
of these equations has a well defined limit in terms of
$\xi$ and $(h_{in},q_{in},Y_{in},Z_{in})$.
This implies that each left-hand side is convergent, uniformly
in $t$ on any compact subset of ${\bf R}$.
Thus, $(h_n,q_n,Y_n,Z_n)$ has a limit
$(h,q,Y,Z)$ in the space
$C^0({\bf R};L^{\infty}_{weak^*}({\bf R};{\bf R}^{1+1+d+d}))$.
This limit satisfies
$$
\int_{-\infty}^{+\infty} {g(s)}{h(t,s)}ds=
\int_{-\infty}^{+\infty} g(\xi(t,y))dy,
$$
$$
\int_{-\infty}^{+\infty}
(q\underline{+}1)(t,s)g(s)ds=
\int_{-\infty}^{+\infty}
(q_{in}\underline{+}1)(s)
g(\xi(t,\xi^{-1}(0,s)\overline{+} t))ds,
$$
$$
\int_{-\infty}^{+\infty}
(Y\overline{+} Z)(t,s)g(s)ds=
\int_{-\infty}^{+\infty}
(Y_{in}\overline{+} Z_{in})(s)
g(\xi(t,\xi^{-1}(0,s)\overline{+} t))ds,
$$
for all functions $g\in L^1({\bf R})$.
\\
Since 
$u=(h,q,Y,Z)$ belongs to
$C^0({\bf R};L^{\infty}_{weak^*}({\bf R};{\bf R}^{1+1+d+d}))$,
we deduce from the previous equations, taken at $t=0$,
that the initial value $u(t=0,\cdot)$
must be equal to $u_{in}=(h_{in},q_{in},Y_{in},Z_{in})$.
Thus, $u=(h,q,Y,Z)$ is a generalized solution to
the augmented system (\ref{non conservative}) in the sense of
(\ref{ww-definition xi},\ref{ww-solution xi},\ref{ww-solution}).
\\
We have proven so far
that, up to a sequence, any sequence $U_n$ in $\Sigma_{\alpha,\delta}$
converges (in the sense of (\ref{convergence}))
to a generalized solution $U=T^{-1}u$.
This solution is valued in 
$CM_{\alpha,\delta}$,  the closed convex hull of $M_{\alpha,\delta}$.
This shows that $\Sigma_{\alpha,\delta}$ is relatively compact
and its closure is contained in the set of generalized solutions
valued in $CM_{\alpha,\delta}$.
\\
Conversely, let us show that all generalized solutions $U$
valued in the closed convex hull 
of $M_{\alpha,\delta}$ belong to the closure of $\Sigma_{\alpha,\delta}$.
Because $T$ preserves convexity, $u=TU$  is
valued in the closed convex hull of $T(M_{\alpha,\delta})$.
Thus, according to a well known property of weak convergence
(see \cite{Ta}, for instance),
the initial value
$u_{in}=u(t=0,\cdot)$ can 
be approached, in the $L^\infty$ weak-* sense,
by a sequence $u_{in,n}$ valued in the manifold $T(M_{\alpha,\delta})$.
Then, we see that the unique generalized solution $U_n$ with initial condition
$T^{-1}u_{in,n}$ must converge  to $U$ (in the sense of 
(\ref{convergence})).
\\
At this point, we have shown that the closure of $\Sigma_{\alpha,\delta}$
is exactly equal to the set of all generalized solutions valued in 
the closed convex hull $CM_{\alpha,\delta}$
of $M_{\alpha,\delta}$.
\\
So, the proof of Theorem \ref{compactness} will be complete
when we are able to show that
$CM_{\alpha,\delta}=CM\cap G_{\alpha,\delta}$.
More concretely, we have to prove that
$$
\{U;\;\;v-\tau\le\alpha-\delta<\alpha+\delta\le v+\tau,
\;\;\;\tau^2+v^2+\eta^2+\zeta^2+2|\tau v-\eta\cdot\zeta|\le 1\}
$$
indeed is the closed convex hull of
$$
\{U;\;\;\delta\le \tau\underline{+}(v-\alpha)\le \frac{1}{\delta}\;;
\;\;\;\tau v=\eta\cdot\zeta,\;\;\tau^2+v^2+\eta^2+\zeta^2=1\}.
$$
We first observe that these sets are equivalently 
defined by
$$
\{U;\;\;\delta\le \tau\underline{+}(v-\alpha)\le \frac{1}{\delta}\;;
\;\;\;(v\underline{+}\tau)^2+|\eta\overline{+}\zeta|^2\le 1\}
$$
and
$$
\{U;\;\;\delta\le \tau\underline{+}(v-\alpha)\le \frac{1}{\delta}\;;
\;\;\;(v\underline{+}\tau)^2+|\eta\overline{+}\zeta|^2= 1\},
$$
respectively. So, the first set, which is compact and convex, certainly
contains the closed convex hull of the second one.
Thus, it is now enough to show that any extremal point $U$ of the first
subset is indeed a point of the second one.
For such a point, for either $\epsilon=1$ or $\epsilon=-1$,
we must have 
$$
(v+\epsilon\tau)^2+|\eta-\epsilon\zeta|^2=1.
$$
Assume $\epsilon=1$ for simplicity, so that
$$
(v+\tau)^2+|\eta-\zeta|^2=1
$$
If $(v-\tau)^2+|\eta+\zeta|^2=1\;\;$, then $U$ belongs to
the second set, as expected. Otherwise, we
have $(v-\tau)^2+|\eta+\zeta|^2<1.$
\\
Let us introduce
$$
U'=(\tau',v',\eta',\zeta')=(0,0, e,-e),
$$
where $e\in {\bf R}^d$, different from zero, is fixed.
we see that for $\lambda\in {\bf R}$ near zero,
the first set still contains $U+\lambda U'$,
which contradicts the assumption
that $U$ is one of its extremal point. 
Indeed, for small $\lambda$, we keep
$$
(v+\lambda v'-\tau-\lambda\tau')^2
+|(\eta+\lambda\eta'+(\zeta+\lambda\zeta')|^2<1,
$$
while we conserve
$$
(v+\lambda v'+\tau+\lambda\tau')^2
+|(\eta+\lambda\eta'-(\zeta-\lambda\zeta')|^2=1,
$$
as well as
$$
\delta\le (\tau+\lambda\tau')\underline{+}
(v+\lambda v'-\alpha)\le \frac{1}{\delta}\;.
$$
The proof of Theorem \ref{compactness}
is now complete.

\subsection*{Acknowledgments}
This article was written at the Bernoulli Centre, 
EPFL, Lausanne, in August 2004,
during the program ``Geometric Mechanics and Its Applications''.
The author is grateful to the organizers, Darryl Holm, Juan-Pablo
Ortega and Tudor Ratiu, for their kind invitation.

This work is also partly supported
by the European IHP project HYKE, HPRN-CT-2002-00282.

\end{document}